\documentclass{ifacconf}

\usepackage{microtype}
\usepackage{graphicx}      
\usepackage{natbib}        
\usepackage{amsmath}       
\usepackage{amssymb}       

\usepackage{mathtools}
\usepackage{graphicx}
\usepackage{physics}

\usepackage{algorithm}
\usepackage{algpseudocode}

\usepackage{spverbatim}

\usepackage{mathtools} 

\usepackage{xcolor}

\usepackage{accents}

\usepackage{subcaption}

\makeatletter
\let\old@ssect\@ssect 
\makeatother

\usepackage[
colorlinks=true
]{hyperref}

\definecolor{dgreen}{rgb}{0, 0.5, 0}

\hypersetup{
colorlinks=true,
linkcolor=blue!50,
citecolor={dgreen},
urlcolor={dgreen},
}

\makeatletter
\def\@ssect#1#2#3#4#5#6{%
  \NR@gettitle{#6}
  \old@ssect{#1}{#2}{#3}{#4}{#5}{#6}
}
\makeatother
\newcommand{\behcet}{Beh\c{c}et}
\newcommand{\acikmese}{A\c{c}\i kme\c{s}e}

\newcommand{\bR}[2]{\mathbb{R}^{#1}_{#2}}
\newcommand{\tcurr}{t_{\mathrm{c}}}
\newcommand{\thorz}{t_{\mathrm{h}}}
\newcommand{\termcost}{L_{\mathrm{h}}}
\newcommand{\termcnstr}{P_{\mathrm{h}}}
\newcommand{\dt}{\mathrm{d}t}
\newcommand{\xicurr}{\xi_{\mathrm{c}}}
\newcommand{\xcurr}{x_{\mathrm{c}}}
\newcommand{\range}[2]{[#1\!:\!#2]}
\newcommand{\tiltermcost}{\tilde{L}_{\mathrm{h}}}
\newcommand{\tiltermcnstr}{\tilde{P}_{\mathrm{h}}}
\newcommand{\wep}{w_{\mathrm{ep}}}
\newcommand{\wprox}{w_{\mathrm{prox}}}

\newcommand{\dtmpc}{\Delta t_{\textsc{mpc}}}

\edef\endfrontmatter{%
  \unexpanded\expandafter{\endfrontmatter}
  \noexpand\endNoHyper 
}
\begin{document}

\begin{frontmatter}
\title{Successive Convexification for\\Nonlinear Model Predictive Control with Continuous-Time Constraint Satisfaction\!\thanksref{footnoteinfo}}

\thanks[footnoteinfo]{\!\!This work was supported by AFOSR grant FA9550-20-1-0053 and ONR grant N00014-20-1-2288. \\ 
\textcopyright 2024 the authors. This work has been accepted to IFAC for publication under a Creative Commons Licence CC-BY-NC-ND.}

\author{$\begin{array}{c}
     \text{\normalsize\textbf{Samet Uzun, Purnanand Elango, Abhinav G.\ Kamath, }} \\ \text{\normalsize\textbf{Taewan Kim, and \behcet~\acikmese}} 
\end{array}$}

\address{William E.\@ Boeing Department of Aeronautics and Astronautics, University of Washington, Seattle, WA 98195, USA \\\normalfont{(e-mail: \texttt{\{samet,\,pelango,\,agkamath,\,twankim,\,behcet\}@uw.edu})}.}
\begin{abstract}
We propose a nonlinear model predictive control (NMPC) framework based on a direct optimal control method that ensures continuous-time constraint satisfaction and accurate evaluation of the running cost, without compromising computational efficiency. We leverage the recently proposed successive convexification framework for trajectory optimization, where: (1) the path constraints and running cost are equivalently reformulated by augmenting the system dynamics, (2) multiple shooting is used for exact discretization, and (3) a convergence-guaranteed sequential convex programming (SCP) algorithm, the prox-linear method, is used to solve the discretized receding-horizon optimal control problems. The 
 resulting NMPC framework is computationally efficient, owing to its support for warm-starting and premature termination of SCP, and its reliance on first-order information only. We demonstrate the effectiveness of the proposed NMPC framework by means of a numerical example with reference-tracking and obstacle avoidance. The implementation is available at \url{https://github.com/UW-ACL/nmpc-ctcs}
\end{abstract}

\end{frontmatter}


\section{Introduction}

\vspace{-0.25cm}
Nonlinear model predictive control (NMPC) is a promising method for automatic control owing to its ability to achieve key mission objectives while satisfying constraints and ensuring robustness to uncertainties \citep{nmpc2000}. This capability mainly stems from leveraging the framework of optimal control, wherein control tasks are specified by means of formulating optimization problems. Then, control laws can be synthesized by solving the formulated problems using appropriate NMPC solvers. Due to these benefits, NMPC has been widely applied to various real-world applications including the control of wind turbines \citep{Schlipf2012}, humanoid robots \citep{tassa2012synthesis}, ground vehicles \citep{di2013stochastic}, and autonomous aerospace systems \citep{eren2017model}.

\vspace{-0.05cm}
Different types of optimal control solvers have been used for NMPC. 
A class of indirect methods based on the maximum principle \citep{liberzon2011calculus} solves a two-point boundary value problem (TPBVP) induced by the necessary conditions for optimality. These indirect methods have been used to apply NMPC to real-world applications \citep{kim2002nonlinear,xie2019pontryagin}. Another class of indirect methods, called differential dynamic programming (DDP), was developed in \citep{mayne1973differential}. The DDP algorithm iteratively performs a forward propagation of the system dynamics and a backward solution update via dynamic programming. The work in \citep{tassa2012synthesis} uses DDP to solve MPC problems with only the system dynamics and the initial condition as constraints. This was later extended to consider input constraints in \citep{tassa2014control} by means of constrained dynamic programming. To handle more general state and input constraints, the augmented Lagrangian method combined with DDP was proposed in \citep{jackson2021altro}.

\vspace{-0.05cm}
The direct optimal control approach to NMPC has also received a lot of attention recently \citep{rawlings2017model}. In this approach, the infinite-dimensional optimal control problem (OCP) at hand is first discretized to yield a finite-dimensional OCP using methods such as single shooting \citep{sargent1978development} or multiple shooting discretization \citep{bock1984multiple, quirynen2015multiple}. The resulting nonlinear problem is then typically solved using nonlinear interior-point methods (nonlinear IPMs) \citep{frison2020hpipm}, sequential quadratic programming (SQP) \citep{houska2011auto}, or sequential convex programming (SCP) \citep{malyuta2022convex}. SQP-type approaches have seen widespread adoption owing to their real-time capability \citep{quirynen2015autogenerating}, especially in conjunction with the real-time iteration (RTI) scheme \citep{diehl2005real, gros2020linear}. 
SCP-based approaches to NMPC have also been gaining popularity recently \citep{dinh2012adjoint, morgan2014model, mao2019convexification}. 

\vspace{-0.05cm}
The SCP approach we adopt has the following advantages over SQP \citep{messerer2021survey, malyuta2021advances}: (1) SCP only requires gradients (first-order information), whereas SQP requires Hessians (second-order information) as well, and (2) SCP subproblems only depend on the primal variables from the previous iteration, whereas SQP subproblems depend on the dual variables as well. 
The NMPC framework we provide consists of $\mathcal{C}^1$-smooth functions, therefore any solution method that requires second-order information such as SQP methods cannot be used; however, any solution method that only requires first-order information can handle the framework we propose.

Typically, NMPC algorithms are limited to nonlinear dynamics with box constraints on the control input \citep{Magni2004}. Further, since piecewise-constant, i.e., zero-order hold, or piecewise-affine, i.e., first-order hold, control input parameterizations are commonly adopted, together with the fact that the control constraints that are levied are convex, these constraints are guaranteed to be satisfied in continuous-time, i.e., intersample constraint satisfaction is guaranteed \citep{szmuk2020successive, kamath2023seco}. However, modern applications require more complicated constraints to be imposed on the control input as well as the state to satisfy mission objectives \citep{Henson1998}. In such situations, the common approach is to only impose constraints at finitely many instants (usually at the node-points within the resulting optimization problem).
However, this approach does not guarantee the satisfaction of constraints in continuous-time. 

A common approach to mitigating this involves using dense grids within the resulting optimization problem or adaptive mesh-refinement techniques \citep{Fontes2019}. However, these approaches lead to increased problem sizes and solve-times, which are in direct conflict with requirements for real-time applications \citep{Zeilinger2014}. Most existing approaches to achieve continuous-time constraint satisfaction involve the computation of maximal robust control invariant sets, mostly in the context of model predictive control with linear systems \citep{Weiss2014, Sadraddini2020} and sampled-data systems \citep{Mitchell2015, Gruber2021}.
Further, a systematic framework to achieve continuous-time constraint satisfaction (CTCS) within model-based NMPC, for real-time applications using online optimization, does not yet exist in the literature, to the best of our knowledge.

To tackle these issues, in this paper, we provide an NMPC framework—that enables continuous-time constraint satisfaction (CTCS) and uses a solution method that can efficiently and reliably solve the formulated problem---by utilizing the real-time-capable successive convexification framework for trajectory optimization that was proposed in \citep{ctcs2024}. Our NMPC framework possesses the following key features:
\begin{enumerate}
    \item A constraint and cost reformulation to enable longer planning horizons with fewer nodes without degradation in inter-sample constraint satisfaction and running cost accuracy.
    \item Multiple shooting exact discretization using the (inverse-free) variational approach \citep{bock1984multiple, quirynen2015autogenerating, kamath2023seco}.
    \item A convergence-guaranteed SCP algorithm, called the prox-linear method \citep{drusvyatskiy2018error}, which enables real-time implementations by supporting warm-starting and premature termination, and due to its reliance on first-order information only.
\end{enumerate}

We showcase the efficacy of the proposed framework through an NMPC simulation that includes reference-tracking and obstacle avoidance.  The implementation is available at
\url{https://github.com/UW-ACL/nmpc-ctcs}
\section{Problem Formulation} \label{sec:probformscp}
\vspace{-0.2cm}
This section presents the successive convexification framework for nonconvex trajectory optimization that was proposed in \cite{ctcs2024}, with a specialization on the receding-horizon optimal control problems encountered in NMPC.
\subsection{Notation}
We adopt the following notation in the subsequent discussion. The set of real numbers is denoted by $\bR{}{}$, the set of nonnegative real numbers by $\bR{}{+}$, and the space of continuous and $k$-times differentiable functions by $\mathcal{C}^k$. The concatenation of vectors $v\in\bR{n}{}$ and $w\in\bR{m}{}$ to provide a vector in $\bR{n+m}{}$ is denoted by $(v,w)$. The set of integers between $a$ and $b$ is denoted by $\range{a}{b}$, where $a\le b$ are integers. The vector of ones in $\bR{n}{}$ is denoted by $1_{n}$, the vector of zeros in $\bR{n}{}$ by $0_n$, the matrix of zeros in $\bR{n\times m}{}$ by $0_{n\times m}$, and the identity matrix in $\bR{n\times n}{}$ by $I_n$. The dimensions are inferred from context wherever omitted. We define $|\square|_+ = \max\{0,\square\}$ and $\|\square\|$ as the 2-norm. Operators $|\square|_+$ and $\square^2$ apply elementwise on a vector.
\subsection{Optimal Control Problem}
We consider a continuous-time dynamical system of the form:
\begin{align}
    \frac{\mathrm{d}\xi(t)}{\dt} = \dot{\xi}(t) = F(t,\xi(t),u(t))
\end{align}
for $t\in\bR{}{+}$, with state $\xi\in\bR{n_{\xi}}{}$ and control input $u\in\bR{n_u}{}$, which is subject to path constraints:
\begin{subequations}
\begin{align}
    g(t,\xi(t),u(t)) \le 0 \\
    h(t,\xi(t),u(t)) = 0
\end{align}\label{path-cnstr}%
\end{subequations}
where $g$ and $h$ are vector-valued functions.

The receding-horizon optimal control problem typical in NMPC can be stated as:
\begin{subequations}
\begin{align}
\underset{u}{\mathrm{minimize}}~~&~\termcost(\xi(\tcurr+\thorz))+\int_{\tcurr}^{\tcurr+\thorz}\!\!\!\!L(t,\xi(t),u(t))\dt\label{ct-cost}\\
\mathrm{subject~to}~&~\dot{\xi}(t) = F(t,\xi(t),u(t)),\quad t\in[\tcurr,\tcurr+\thorz]\\
~&~g(t,\xi(t),u(t)) \le 0,~\:\!\qquad t\in[\tcurr,\tcurr+\thorz]\\
~&~h(t,\xi(t),u(t)) = 0,~\qquad t\in[\tcurr,\tcurr+\thorz]\\
~&~\xi(\tcurr) = \xicurr\\
~&~\termcnstr(\xi(\tcurr+\thorz)) \le 0
\end{align}\label{ct-ocp}%
\end{subequations}
where $\tcurr\in\bR{}{+}$ is the current time, $\thorz\in\bR{}{+}$ is the planning horizon time duration, $L$ is the running cost function, $\termcost$ is the terminal cost function, $\xicurr$ is the current state, and $\termcnstr$ is the terminal constraint function. We assume that functions $F$, $g$, $h$, $L$, $\termcost$, and $\termcnstr$ are $\mathcal{C}^1$ (and potentially nonconvex).
\subsection{Reformulation of Constraints and Running Cost}
We convert \eqref{ct-ocp} to the Mayer form \citep[Sec. 3.3.2]{liberzon2011calculus} by reformulating the path constraints and the running cost as auxiliary dynamical systems—as follows:
\begin{align}
    \dot{x}(t) & = \left[ \begin{array}{c} \dot{\xi}(t) \\[0.05cm] \dot{l}(t) \\[0.05cm] \dot{y}(t) \end{array} \right]\nonumber\\
    & = \left[ \begin{array}{c} F(t,\xi(t),u(t)) \\[0.05cm] L(t,\xi(t),u(t)) \\[0.05cm] 1^\top |g(t,\xi(t),u(t))|_+^2 + 1^\top h(t,\xi(t),u(t))^2 \end{array} \right]\nonumber\\  
    & = f(t,x(t),u(t))\label{ct-reform-dyn}
\end{align}
The integral of the running cost function in \eqref{ct-cost} is equal to $l(\tcurr+\thorz)$ when the initial condition $l(\tcurr) = 0$ is imposed. In addition, a periodic boundary condition on the constraint violation integrator:
\begin{align}
    y(\tcurr) = y(\tcurr+\thorz)\label{cnstr-bc}
\end{align}
is specified to ensure that the path constraints are satisfied at all times. 

Furthermore, we define:
\begin{subequations}
\begin{align}
    \tiltermcost(x(\tcurr+\thorz)) ={} & \termcost(\xi(\tcurr+\thorz)) + l(\tcurr+\thorz)\\
    \tiltermcnstr(x(\tcurr+\thorz)) ={} & \termcnstr(\xi(\tcurr+\thorz))
\end{align}
\end{subequations}
\subsection{Discretization} \label{sec:disc}
Discretization is a crucial step in direct numerical optimal control for computational tractability. We adopt the multiple-shooting approach \citep{bock1984multiple, quirynen2015autogenerating, kamath2023seco} for discretizing \eqref{ct-reform-dyn}. The first step is to generate a uniform grid of $N$ nodes within the interval $[\tcurr,\tcurr+\thorz]$:
$$
    \tcurr = t_1 < \ldots < t_N = \tcurr + \thorz
$$
The interval $[t_k,t_{k+1}]$, for $k\in\range{1}{N\!-\!1}$, is called a sub-interval, and its length is denoted as: 
\begin{equation} \label{eq:sub_int}
    \Delta t = \thorz/(N-1)
\end{equation}
Next, we parameterize the control input with a first-order-hold (FOH):
\begin{align}
    u(t) = \left( \frac{t_{k+1} - t}{t_{k+1}-t_k} \right) u_k + \left( \frac{t - t_{k}}{t_{k+1}-t_k} \right)u_{k+1}\label{FOH-param}
\end{align}
for $t\in[\tcurr,\tcurr+\thorz]$ and $k\in\range{1}{N\!-\!1}$. Given the choice of parameterization \eqref{FOH-param}, convex constraints on the control input encoded in \eqref{path-cnstr} can be separated out and represented with a convex set $\mathcal{U}$. Continuous-time satisfaction of the convex control constraints is ensured by constraining $u_k$ alone, i.e., $u_k\in\mathcal{U}$, for $k\in\range{1}{N}$. In practice, the combination of time-dilation and FOH parameterization of the control input is sufficient for accurately modeling a large class of optimal control problems. See \cite{malyuta2019discretization} for a detailed discussion on the relative merits of popular discretization and parameterization techniques in numerical optimal control.

We treat the node-point values of the state and control input, $x_k$ and $u_k$, respectively, as decision variables, and re-write the system dynamics \eqref{ct-reform-dyn} equivalently as:
\begin{align}
    x_{k+1} = x_{k} + \int_{t_k}^{t_{k+1}} f(t,x^k(t),u(t))\dt \label{dt-dyn}
\end{align}
for $k\in\range{1}{N\!-\!1}$, where $x^k(t)$ is the solution to \eqref{ct-reform-dyn} over $[t_k,t_{k+1}]$ with initial condition $x_k$, and $u(t)$ is FOH-parameterized using $u_k$ and $u_{k+1}$. Note that the mapping $t\mapsto f(t,x^k(t),u(t))$ is piecewise-linear, hence it is Riemann integrable. We use a variety of off-the-shelf ODE solvers to compute the integral in \eqref{dt-dyn} to at most machine precision. 

Then, \eqref{dt-dyn} can be compactly expressed 
as follows:
\begin{align}
    x_{k+1} = f_k(x_k,u_k,u_{k+1})
\end{align}
The boundary condition \eqref{cnstr-bc} is relaxed to:
\begin{align} \label{eq:licq}
    y_{k+1} - y_{k} \le \epsilon
\end{align}
for $k\in\range{1}{N\!-\!1}$ and a sufficiently small $\epsilon>0$, so that all feasible solutions do not violate the linear independence constraint qualification (LICQ) \citep{ctcs2024}. Note that $y_k$ is the value of the state $y$ in \eqref{ct-reform-dyn} at node $t_k$.

\vspace{-0.05cm}
The discretized form of receding-horizon optimal control problem \eqref{ct-ocp}, with the reformulation of constraints and running cost, and with the parameterized control input, can be stated as:
\begin{subequations}
\begin{align}
\underset{x_k,u_k}{\mathrm{minimize}}~~&~\tilde{L}_{\mathrm{h}}(x_N) \label{dt-ocp:cost} \\
\mathrm{subject~to}~&~x_{k+1} = f_k(x_k,u_k,u_{k+1}),~~ k \in \range{1}{N\!-\!1} \label{dt-ocp:dyn}\\
~&~E_y(x_{k+1}-x_k) \le \epsilon, ~~\,\qquad k \in \range{1}{N\!-\!1} \label{dt-ocp:relax}\\
~&~u_k\in\mathcal{U}, \qquad\qquad\qquad\qquad k \in\range{1}{N}\label{dt-ocp:ctrl}\\
~&~x_1 = \xcurr\label{dt-ocp:initbc}\\
~&~\tiltermcnstr(x_N) \le 0\label{dt-ocp:termbc}
\end{align}\label{dt-ocp}%
\end{subequations}
where $\xcurr = (\xicurr,0,0)$ and $E_y = \big[0_{n_{\xi}}^\top~0~1\big]$.
\vspace{-0.15cm}
\subsection{Sequential Convex Programming} \label{subsec:scp}
\vspace{-0.25cm}
We adopt the prox-linear method \citep{drusvyatskiy2018error}—which is closely related to the penalized trust region (PTR) algorithm for trajectory optimization \citep{szmuk2020successive,reynolds2020real}—to solve:
\begin{subequations}
\begin{align}
\underset{x_k,u_k}{\mathrm{minimize}}~~&~\tilde{L}_{\mathrm{h}}(x_N)\label{dt-penal:cost}\\[-0.25cm]
 &+\wep\sum_{k=1}^{N-1}\|x_{k+1} - f_k(x_k,u_k,u_{k+1})\|_1  \nonumber\\[-0.1cm]
 &+\wep\max\big\{0,\tiltermcnstr(x_N)\big\}\nonumber\\
\mathrm{subject~to}~&~E_y(x_{k+1}-x_k) \le \epsilon, ~~\,\qquad k \in \range{1}{N\!-\!1}\label{dt-penal:relax}\\
~&~u_k\in\mathcal{U}, \qquad\qquad\qquad\qquad k \in\range{1}{N}\label{dt-penal:ctrl}\\
~&~x_1 = \xcurr\label{dt-penal:initbc}
\end{align}\label{dt-penal-ocp}%
\end{subequations}
where the nonconvex constraints in \eqref{dt-ocp} are exactly penalized with the $\ell_1$-penalization. For a large enough, finite $\wep$, we can show that a solution of \eqref{dt-penal-ocp} that is feasible with respect to \eqref{dt-ocp:dyn} and \eqref{dt-ocp:termbc} is also a solution to \eqref{dt-ocp}. See \cite[Chap. 17]{nocedal2006numerical} for a detailed discussion.

\vspace{-0.05cm}
The penalized problem \eqref{dt-penal-ocp} can be compactly stated as:
\begin{align}
    \underset{Z\in\mathcal{Z}}{\mathrm{minimize}}~~G(Z) + \wep H(C(Z))\label{dt-penal-ocp-compact}
\end{align}
where $Z = (x_1,\ldots,x_N,u_1,\ldots,u_N)$, $\mathcal{Z}$ is a convex set corresponding to constraints \eqref{dt-penal:relax}, \eqref{dt-penal:ctrl}, and \eqref{dt-penal:initbc}, $G$ is a proper convex function, and 
$H \circ C$ is a composition of a proper convex function and $\mathcal{C}^1$-smooth mapping corresponding to \eqref{dt-penal:cost}.
In the case where 
$\tilde{L}_{\mathrm{h}}$ is convex, it is incorporated into $G$, while in the case where $\tiltermcnstr$ is convex, constraint \eqref{dt-ocp:termbc} is integrated into $\mathcal{Z}$ without any penalization.
It is also possible to specify a different exact penalty weight for each of the penalized constraints.

\vspace{-0.05cm}
The prox-linear method determines a stationary point of \eqref{dt-penal-ocp-compact} by solving a sequence of convex subproblems. At iteration $j+1$, it solves the following convex problem:
\begin{align}
    \underset{Z\in\mathcal{Z}}{\mathrm{minimize}}~&G(Z) + \wep H(C(Z^{j})+\nabla C(Z^{j})(Z-Z^{j}))\label{cvx-subproblem}\\[-0.15cm]
    & + \frac{\wprox}{2}\|Z-Z^j\|_2^2\nonumber
\end{align}
where $Z^j$ is the solution from iteration $j$. The sequence $Z^j$ can be shown to converge for an appropriate choice of the proximal term weight $\wprox$. Further, if the converged point is feasible for \eqref{dt-ocp}, then it is also a KKT point. The user-specified initialization to the algorithm is $Z^1$. In the subsequent section, we will discuss a warm-start strategy that provides an approximate solution as initialization. 

The specialization of \eqref{cvx-subproblem} at iteration $j+1$ to \eqref{dt-penal-ocp} can be stated as:
\begin{subequations}
\begin{align}
    \underset{x_k,u_k}{\mathrm{minimize}}~~&~\tiltermcost(x^j_N) + \nabla\tiltermcost(x^j_N)(x_N-x^j_N)\label{cvx-subproblem-custom:obj}\\[-0.3cm]
    &+ \wep\nu_{\mathrm{h}} + \wep\sum_{k=1}^{N-1}1_{n_x}^\top(\mu^+_k+\mu_k^-) \nonumber \\
    &+ \frac{\wprox}{2}\sum_{k=1}^N \|x_k-x^j_k\|_2^2 + \|u_k - u^j_k\|_2^2 \nonumber \\[-0.05cm]
    \mathrm{subject~to}~&~x_{k+1} = A_kx_k + B^-_ku_k + B^+_ku_{k+1}\label{cvx-subproblem-custom:dyn}\\
     &\,\quad\qquad+ \mu_k^+ - \mu_k^- + w_k,~~ k\in\range{1}{N\!-\!1}\nonumber\\
     &~E_{y}(x_{k+1}-x_{k}) \le \epsilon, ~~\,\qquad k\in\range{1}{N\!-\!1}\label{cvx-subproblem-custom:relax}\\
     &~\mu_k^+\ge 0,~\mu^-_k\ge 0 \!\,~\qquad\qquad k\in\range{1}{N\!-\!1} \\
     &~u_k \in \mathcal{U} \!\,~\qquad\qquad\qquad\qquad k\in\range{1}{N}\\
     &~x_1 = \xcurr\\
     &~\tiltermcnstr(x^j_N) + \nabla\tiltermcnstr(x^j_N)(x_N-x^j_N) \le \nu_{\mathrm{h}}\\
     &~\nu_{\mathrm{h}} \ge 0
\end{align}\label{cvx-subproblem-custom}%
\end{subequations}
where slack variables $\mu_k^+$ and $\mu^-_k$ are introduced through an equivalent reformulation of the $\|\square\|_1$ term, and $\nu_{\mathrm{h}}$ through the reformulation of the $\max\{0,\square\}$ term in the objective function \eqref{cvx-subproblem-custom:obj}. Furthermore, $A_k$, $B^-_k$, $B^+_k$, and $w_k$, for $k\in\range{1}{N\!-\!1}$, form the gradient of \eqref{dt-ocp:dyn} at $Z^j = (x^j_1,\ldots,x^j_N,u^j_1,\ldots,u^j_N)$. More precisely, for each $k\in\range{1}{N\!-\!1}$, let $x^{jk}(t)$ denote the solution to \eqref{dt-dyn} over $[t_k,t_{k+1}]$ generated with initial condition $x^j_k$ and control input $u^j(t)$, which is FOH-parameterized using $u^j_k$ and $u^j_{k+1}$. The Jacobians of $f$, evaluated with respect to $x^{jk}(t)$ and $u^j(t)$, are denoted by:
\begin{align*}
    A^k(t) ={} & \left.\frac{\partial f(x,u)}{\partial x}\right|_{(x^{jk}(t),u^j(t))}\\
    B^k(t) ={} & \left.\frac{\partial f(x,u)}{\partial u}\right|_{(x^{jk}(t),u^j(t))}
\end{align*}
Next, we solve the following initial value problem:
\begin{align*}
    \dot{\Phi}_x(t,t_k) ={} & A^k(t)\Phi_x(t,t_k)\\
    \dot{\Phi}^{-}_u(t,t_k) ={} & A^k(t)\Phi^{-}(t,t_k) + B^k(t)\left(\frac{t_{k+1}-t}{t_{k+1}-t_k}\right)\\
    \dot{\Phi}^{+}_u(t,t_k) ={} & A^k(t)\Phi^{+}(t,t_k) + B^k(t)\left(\frac{t-t_k}{t_{k+1}-t_k}\right)\\
    \Phi_x(t_k,t_k) ={} & I_{n_x}\\
    \Phi_u^{-}(t_k,t_k) ={} & 0_{n_x\times n_u}\\
    \Phi_u^{+}(t_k,t_k) ={} & 0_{n_x\times n_u}
\end{align*}
to obtain:
\begin{align*}
    A_k ={} & \Phi_x(t_{k+1},t_k)\\
    B_k^{-} ={} & \Phi_u^{-}(t_{k+1},t_k)\\
    B_k^{+} ={} & \Phi_u^{+}(t_{k+1},t_k)\\
    w_k ={} & x^{jk}(t_{k+1}) - A_k x^j_k - B^-_k u^j_k - B_k^{+}u^{j}_{k+1} 
\end{align*}

See \cite[Sec. C]{chari2024fast} for further details about gradient computation, and numerical scaling and preconditioning operations that are beneficial in terms of obtaining reliable numerical performance of SCP in practice.



We use CVXPY \citep{diamond2016cvxpy} with either ECOS \citep{domahidi2013ecos} or MOSEK \citep{aps2019mosek} for modeling and solving \eqref{cvx-subproblem-custom}. Customizable conic solvers such as the proportional-integral projected gradient (PIPG) algorithm \citep{yu2022proportional, yu2022extrapolated} could be used for efficient implementations.
Observe that the structure of \eqref{cvx-subproblem-custom} is unaffected by the type and number of path constraints \eqref{path-cnstr}. As a result, the effort of parsing \eqref{cvx-subproblem-custom} to the canonical form of a convex optimization solver is the same for a large class of optimal control problems. This is particularly beneficial for producing efficient customized implementations \citep{Elango2022, Kamath2023} for solving \eqref{cvx-subproblem-custom}, which will apply to a variety of scenarios without the need for any additional modifications. Further, this allows for enabling or disabling constraints within the OCP without the need for re-parsing or re-customization, since the subproblem structure remains intact.

We define the following quantities to later assess the performance of SCP within a NMPC simulation:
\begin{subequations} \label{eq:scp_performance}
\begin{align}
    J_l ={} & E_l\sum_{k=1}^{N-1}\mu_k^++\mu_k^-\\
    J_y ={} & E_y\sum_{k=1}^{N-1}\mu_k^++\mu_k^-\\
    J_{\mathrm{prox}} ={} & \sum_{k=1}^N \|x_k-x_k^j\|_2^2 + \|u_k-u_k^j\|_2^2
\end{align}    
\end{subequations}
where $E_l = \big[0_{n_{\xi}}^\top~1~0\big]$. Taken together, these quantities measure the linearization error in computing the integral of the running cost, violation of the path constraints, and overall convergence of the SCP.


\section{NMPC Framework}

This section outlines the nonlinear model predictive control (NMPC) framework adopted in this paper, elucidating the specific design decisions made.
NMPC recursively solves the receding-horizon optimal control problem \eqref{ct-ocp}. This method is applied in scenarios characterized by dynamic conditions, encompassing uncertain states, evolving objectives, and varying constraints. The objective and constraints involved might exhibit state-dependent or time-varying characteristics, with uncertain temporal patterns. NMPC offers a certain level of robustness against these uncertainties and varying characteristics by recursively solving the receding-horizon optimal control problem \citep{nmpc2000}.

As described in Section \ref{sec:probformscp}, we first reformulate the constraints and running cost in \eqref{ct-ocp}, parameterize the infinite-dimensional control input, discretize the continuous-time state for numerical tractability, and then employ an SCP algorithm to solve the discretized optimal control problem.    
The reformulation process yields an augmented dynamical system where the right-hand-side function is $\mathcal{C}^1$-smooth, and the primary advantage of the SCP framework originates from its reliance solely on first-order information.
Reformulating constraints and the running cost to embed them within the dynamics, as shown in \eqref{ct-reform-dyn}, offers the benefit of being able to operate on sparse discretization grids without degradation in inter-sample constraint satisfaction and running cost accuracy.
This approach equips the NMPC framework with the capability to handle longer horizons while using fewer discretization nodes.

\begin{figure}[!htpb]
\centering
\begin{subfigure}[b]{0.99\linewidth}
\includegraphics[width=\linewidth]{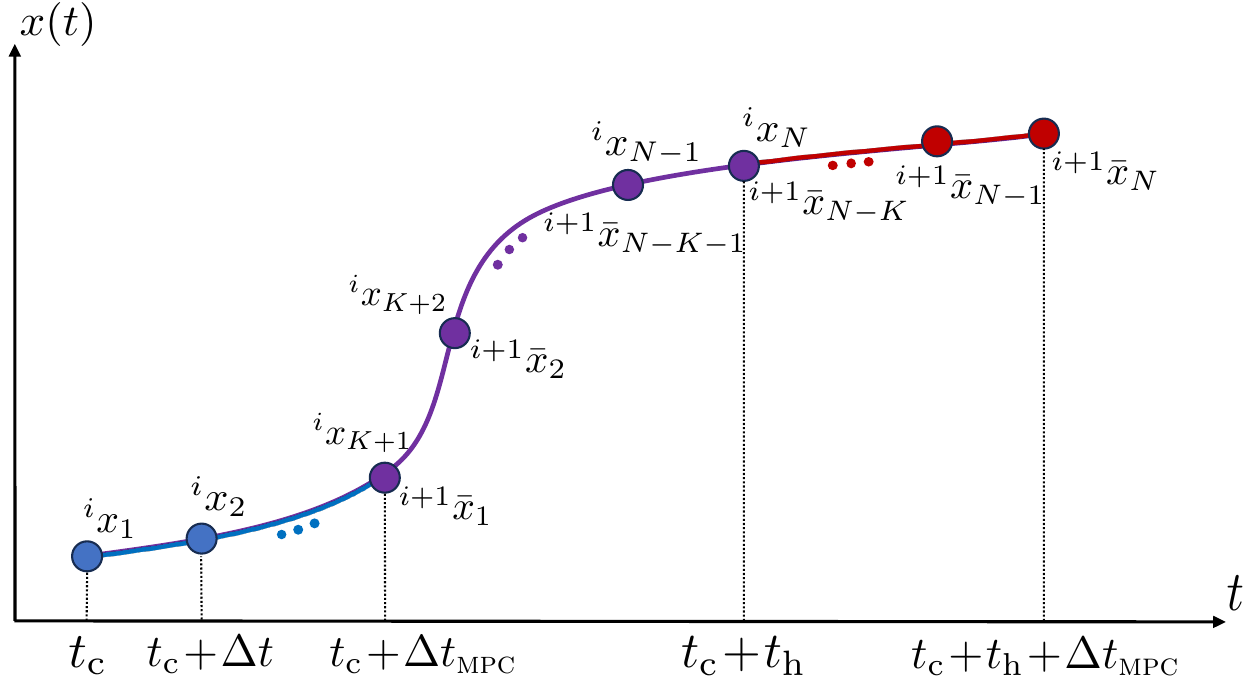}
\end{subfigure}
\caption{An illustration of the warm-starting procedure. The blue and purple trajectories form the SCP solution at $i$th NMPC run. The purple and red trajectories form the initialization to SCP at $(i+1)$th NMPC run.}
\label{fig:nmpc_init}
\end{figure}

A crucial parameter to consider in addition to the number of node points is the time interval between each NMPC run, denoted as $\dtmpc = K \Delta t$, where $\Delta t$ is the discretization grid sub-interval length \eqref{eq:sub_int}, and $K$ is a positive real number satisfying $K \leq N-1$.
When determining suitable values for parameters $N$ and $K$, it is crucial to consider the degrees of freedom required in the parameterized control input and the stiffness \citep{wanner1996solving} of the dynamical system \eqref{ct-reform-dyn}.
Given that the degree of freedom of the parameterized control input is directly related to $N$, it should be selected to be sufficiently large to ensure the feasibility of the discretized optimal control problem \eqref{dt-ocp}.
In situations where the system's state undergoes rapid changes or when objectives and constraints rapidly vary, the dynamical system \eqref{ct-reform-dyn} tends to become stiff. As a consequence:
\begin{enumerate}
    \item Selecting $N$ to yield a suitably small $\Delta t$ is desirable to alleviate challenges related to ill-conditioned sensitivity matrices in stiff systems. Additionally, opting for a small $\dtmpc$ enables NMPC to respond to changes promptly. Hence, $K$ cannot be a large value.
    \item In non-stiff systems, the reformulation of the optimal control problem enables the selection of the number of node points based on the number of degrees of freedom required in the parameterized control input without risking inter-sample constraint violations or compromising the accuracy of the running cost.
    Furthermore, for computational efficiency, $K$ can be increased as long as $\dtmpc$ is not larger than the fastest time-scales of the system.
\end{enumerate}
Another important step for computational efficiency is warm-starting. Given that NMPC recursively solves \eqref{dt-ocp}, an effective warm-starting technique that can be applied after the first run of the NMPC, can allow SCP to converge to an acceptable solution within a limited number of iterations.
The left superscript of the variables indicates the specific MPC run, while the right superscript denotes the SCP iteration described in Section \ref{sec:probformscp}.
When $K$ is a positive integer, SCP at the current time instant is warm-started by shifting the SCP solution from the previous time instant. An example illustrating the warm-starting procedure is shown in Figure \ref{fig:nmpc_init}.

The warm-start solution uses $N-K$ nodes of the preceding SCP solution as the initialization for the first $N-K$ nodes of the subsequent SCP call:
\vspace{0.08cm}
\begin{align*}
    {}^{i+1}\bar{x}_j &= {}^i x_{j+K}\\
    {}^{i+1}\bar{u}_j &= {}^i u_{j+K}
\end{align*}
for $j \in \range{1}{N-K}$. The initialization for the control input at the last $K$ nodes is set to the control input at the $(N-K)$th node of the previous SCP solution:
\vspace{0.08cm}
\begin{align*}
    {}^{i+1}\bar{u}_{N-K+j} &= {}^i u_{N}
\end{align*}
for $j \in \range{1}{K}$. The initialization for the states at the last $K$ nodes is determined by integrating dynamics \eqref{ct-reform-dyn} from the last state of the previous SCP solution with the control inputs at the last $K$ nodes as follows:
\vspace{0.05cm}
\begin{align*}
   \!\! & {}^{i+1}\bar{x}_{N-K+j}  {}= {}^{i}x_{N} + \int_{{\tcurr}+\thorz}^{\tcurr+\thorz+j \Delta t}\!\!\!\!\!f(t, x(t), {}^i u_N) \dt
\end{align*}
for $j \in \range{0}{K-1}$. Therefore, if the preceding SCP solution is dynamically feasible, then the initialization for the subsequent SCP call will be dynamically feasible as well. This aspect of warm-starting is crucial for achieving rapid convergence.
Due to this rapid convergence facilitated by warm-starting, we terminate SCP prior to complete convergence (similar to the real-time-iteration (RTI) scheme in \citep{diehl2005real}). We confirm through numerical experiments that premature termination does not compromise the quality of the solutions; however, it notably enhances computational efficiency.

The overall NMPC framework is summarized in Algorithm \ref{alg:NMPC}. 
The inputs to Algorithm \ref{alg:NMPC} comprise the total execution time of NMPC, $T$, the horizon time $\thorz$, and grid sizes $N$ and $K$ for defining $\Delta t$ and $\Delta t_{\textsc{mpc}}$, along with an initialization for the first SCP call denoted as:
\vspace{0.05cm}
$$
    {}^1\!\bar{Z} =  ({}^1\bar{x}_1,\ldots,{}^1\bar{x}_N,{}^1\bar{u}_1,\ldots,{}^1\bar{u}_N)
$$
Given the maximum number of iterations $j_{\max}$, grid size $N$, current time $\tcurr$, and time duration $\thorz$, we define function $\textsc{scp}_{N,\tcurr,\thorz}^{j_{\max}}$, which applies the SCP algorithm described in Section \ref{subsec:scp} to solve \eqref{dt-ocp} with inputs: current state estimate $\hat{x}_{\mathrm{c}}$, running cost function $L$, and path constraint functions $g$ and $h$. Note that the running cost function captures the tracking cost for the current segment of a reference trajectory over $[\tcurr,\tcurr+\thorz]$.
Each run of NMPC starts with information acquisition by the sensors, which then provides estimates of the current state, running cost, and constraints, all of which dynamically evolve over time.
Then, the SCP call yields a sub-optimal solution for \eqref{dt-ocp}, which is also used to warm-start SCP in the subsequent NMPC run.
The control input solution from the SCP call is applied to the system until the time instant of the next NMPC run.

%
\begin{algorithm}[!htpb] 
\caption{NMPC Framework}
\label{alg:NMPC}
\begin{algorithmic}[0] 
\State \textbf{Inputs:} 
$T$, $\thorz$, $K$, $N$, $j_{\max}$, ${}^{1}\!\bar{Z}$\\
$i \gets 0$\\
$\tcurr \gets 0$\\
$\Delta t \gets \thorz / (N-1)$\\
$\Delta t_{\textsc{mpc}} \gets K \Delta t$
\While{$ \tcurr \leq T $}
\State $\hat{x}(\tcurr), {}^{i}\!L, {}^{i}\!g, {}^{i}\!h \gets$ \textsc{sensor}($x(\tcurr)$)
\State ${}^{i} Z \gets \textsc{scp}_{N,\tcurr,\thorz}^{j_{\max}}(\hat{x}(\tcurr), {}^{i}\!L, {}^{i}\!g, {}^{i}\!h, {}^{i}\!\bar{Z})$
\State ${}^{i+1} \bar{Z} \gets$ \textsc{warm\_start}$({}^iZ)$
\State $x(\tcurr + \dtmpc) \gets $ \textsc{evolve}($x(\tcurr)$, ${}^{i}\!Z$)
\State $i \gets i+1$
\State $\tcurr \gets \tcurr + \dtmpc$
\EndWhile 
\end{algorithmic}
\end{algorithm}
\section{Numerical results}

\vspace{-0.2cm}
This section describes a numerical demonstration of the proposed NMPC framework for output reference tracking, with obstacle avoidance, for a double-integrator with drag and bounded control inputs. The system state is defined as $\xi(t) = (r(t),v(t))$, where $r(t)$ and $v(t)$ are the vehicle's position and velocity, and the right-hand-side of the governing ODE is given by:
\vspace{0.2cm}
$$ F(t,\xi(t),u(t)) =     
\begin{bmatrix}
    v(t) \\
    u(t) - c_{\mathrm{d}} \| v(t) \| v(t)
\end{bmatrix}
$$
where $c_{\mathrm{d}}$ is the drag coefficient of the vehicle.

The running cost is defined as:
\vspace{0.1cm}
$$ L(t,\xi(t),u(t)) = \|r(t) - r_{\mathrm{ref}}(t)\|^2$$
where $r_{\mathrm{ref}}$ is the reference output.
We assume that the vehicle has a sensor, such as a LiDAR, that detects the obstacle within a range denoted by $r_{\mathrm{sens}} \in \mathbb{R}_+$. We include the instantaneously detected obstacles to \eqref{ct-reform-dyn} at each NMPC run.
We assume that the obstacles are elliptical and define an index set of the obstacles that are detected by the sensor at time $\tcurr$ as:
\vspace{0.2cm}
$$
\mathcal{O}_{\mathrm{c}} = \{ i \mid \| S^i ( r(\tcurr) - r^i_{\mathrm{obs}} )\| \leq r_{\mathrm{sens}},~i \in \range{1}{n_{\mathrm{obs}}}\}
$$
where 
$n_{\mathrm{obs}}$ is the total number of the obstacles, $S^i$ is the shape matrix of $i$th obstacle, $r_{\mathrm{obs}}$ is a matrix such that the $i$th column of the matrix, $r^i_{\mathrm{obs}}$, represents the center of the $i$th obstacle.
The vector-valued path constraint function, which includes the obstacle avoidance constraints, is given by:
\vspace{0.2cm}
$$
g(t,\xi(t),u(t)) = 
    \begin{bmatrix}
        \vdots \\
        1-\|S^i(r(t) - r^i_{\mathrm{obs}})\| \\
        \vdots
    \end{bmatrix},~\forall i \in \mathcal{O}_\mathrm{c}
$$
The boundary condition defined for CTCS in (\ref{cnstr-bc}) is relaxed by $\epsilon > 0$ in (\ref{eq:licq}) to satisfy LICQ. However, this relaxation may cause a pointwise violation of the constraint. In \citep{ctcs2024}, it is proved that for each relaxation parameter $\epsilon > 0$, there exists a constraint tightening value $\delta_{\textrm{obs}} > 0$ such that 
\vspace{0.2cm}
\begin{align*}
& \int_{t_k}^{t_{k+1}} | g_i(t, \xi(t), u(t))+\delta_{\textrm{obs}}|_{+}^2\mathrm{d}t \le \epsilon \\
& \iff \int_{t_k}^{t_{k+1}} |g_i(t, \xi(t), u(t))|_{+}^2\mathrm{d}t = 0
\end{align*}
The choice of constraint tightening value $\delta_{\textrm{obs}}$ is determined by the results of numerical experiments.

Note that the obstacle avoidance constraint function is not differentiable at $0$, specifically when the vehicle's position coincides with the center point of the obstacle.
However, this situation is not typically expected to occur unless the initialization is intentionally configured this way.
If this situation presents an issue, we can address it by using a smooth quadratic in place of the norm.

Convex constraints on the control input are encoded by $\mathcal{U} = \{ u\in\bR{2}{}\,|\, \|u\| \le u_{\max} \}$.
Table \ref{tab:dyn-obstacle-avoidance-param} shows the system and simulation parameters.
\begin{table}[!htpb]
\centering
\caption{}\label{tab:dyn-obstacle-avoidance-param}
{\renewcommand{\arraystretch}{1.1}
\begin{tabular}{l|l}
\hline
Parameter & Value\\\hline\\[-0.3cm]
$T$, $t_{\mathrm{h}}$, $K$, $N$, $j_{\mathrm{max}}$
& $30$ s, $8$ s, $4$, $17$, $3$ \\
$u_{\max}$ & $5$ m s$^{-2}$ \\
$c_{\mathrm{d}}$ & $0.05$ m$^{-1}$ \\
$r_{\mathrm{ref}}(t)$ & 
$\left[ \begin{array}{c}
-10 + 2t/3 \\
5 e^{-0.05t} \sin \left( \pi/2 + t/5 + t^2/36 \right)
\end{array} \right]$\\[0.4cm]
$S^i$, $i=1,\dots,5$ & 
$\left[ \begin{array}{cc}
2.4 & 0\\
0   & 0.4
\end{array} \right]$ \\[0.4cm]
$r_{\mathrm{obs}}$ &
$\left[ \begin{array}{cccccccccc}
       -6.5 & -1.8 & 0.5 & 2.7 & 6.2 \\
       -3 & 3.5 & -3.5 & 1.5 & -3
    \end{array} \right]$ \\[0.4cm]
$\delta_{\textrm{obs}}$ & $0.05$ \\
$r_{\mathrm{sens}}$ & $6$ m \\ [0.1cm]
\hline
\end{tabular}}
\end{table}

We study the impact of premature termination on the quality of the solution.
Within each NMPC run, we run SCP for $10$ iterations, and in each SCP iteration, we record the values of the $J_l$, $J_y$, and $J_{\mathrm{prox}}$ terms defined in \eqref{eq:scp_performance}.
When these quantities converge to the chosen convergence tolerance, it indicates that the integral of the running cost in the reformulated problem is accurate, the obstacle avoidance constraint is satisfied, and the overall convergence of SCP.

\begin{figure}[!htpb]
\centering
\begin{subfigure}[b]{0.99\linewidth}
\includegraphics[width=\linewidth]{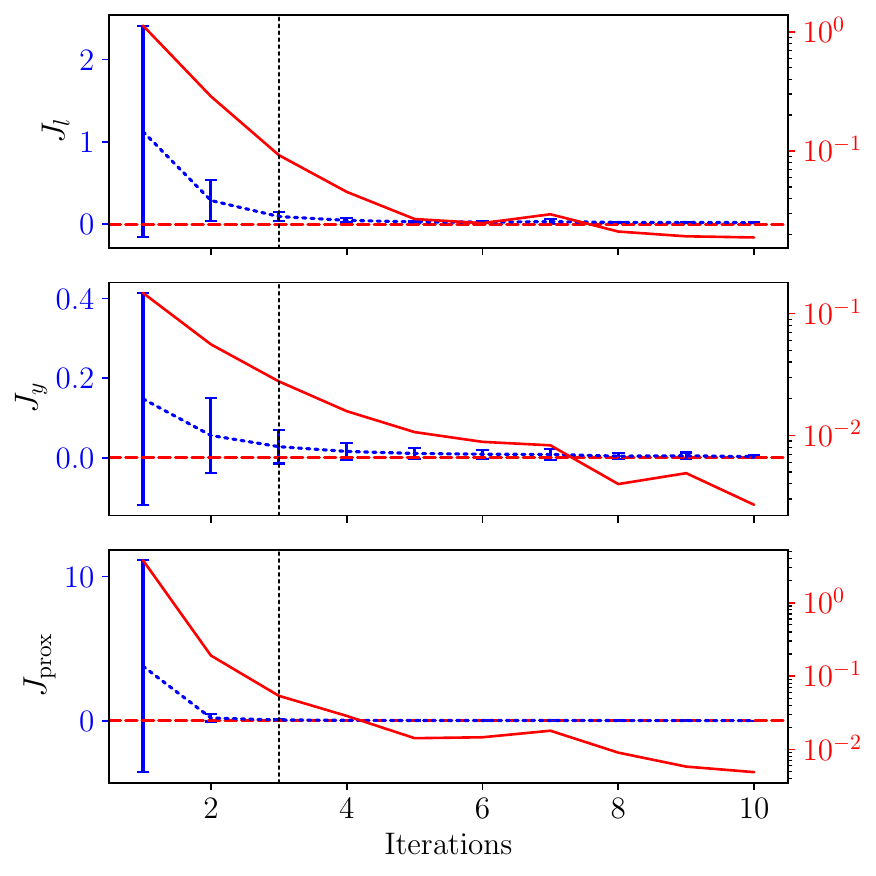}
\end{subfigure}
\caption{The mean and standard deviation of $J_l$, $J_y$, and $J_{\mathrm{prox}}$ along the trajectory for $10$ SCP iterations. The linear scale on the left axis corresponds to the mean and standard deviation shown by the blue curve, while the log scale on the right axis corresponds to the mean shown by the red curve. The black-dashed line denotes the iteration at which the SCP algorithm is set to terminate.
}
\label{fig:simerrors}
\end{figure}

Figure \ref{fig:simerrors} shows the variation in the mean and standard deviation of the values of these quantities along the trajectory as a function of the number of SCP iterations.
These plots indicate that terminating SCP after $3$ iterations yields an acceptable solution.

\begin{figure}[!htpb]
\centering
\begin{subfigure}[b]{0.99\linewidth}
\includegraphics[width=\linewidth]{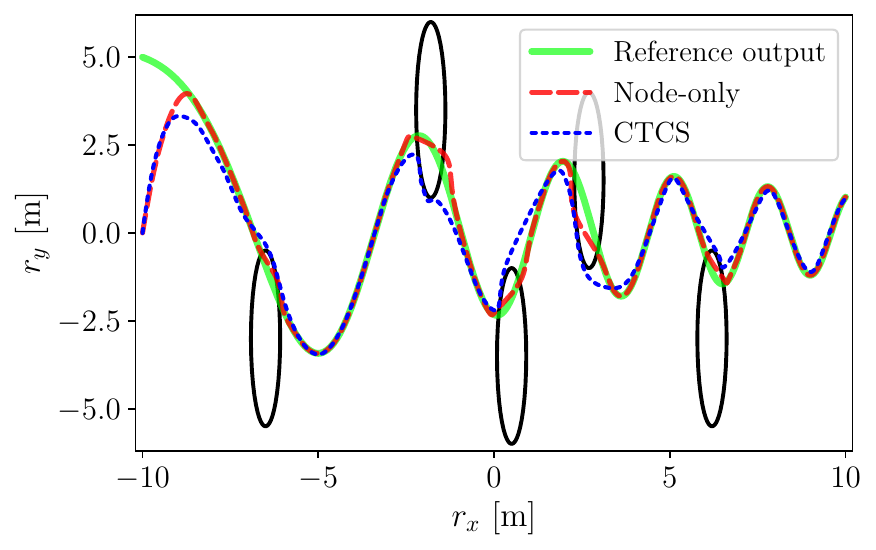}
\end{subfigure}
\caption{$r_x$ and $r_y$ represent the first and second components of the position vector $r$, respectively. A green solid curve illustrates the reference output, while black ellipses represent the obstacles. The trajectories obtained from NMPC using the node-only and the CTCS formulations are presented as red-dashed and blue-dotted curves, respectively. 
}
\label{fig:sim_path_thick}
\end{figure}

\begin{figure}[!htpb]
\centering
\begin{subfigure}[b]{0.99\linewidth}
\includegraphics[width=\linewidth]{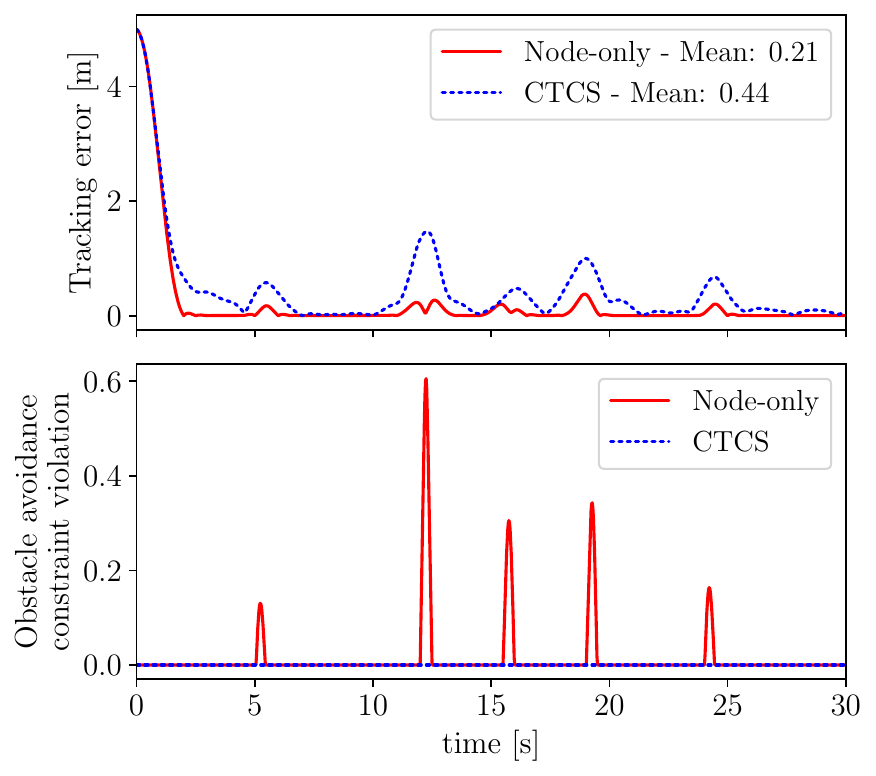}
\end{subfigure}
\caption{The first plot shows the discrepancy between the reference output and the trajectories obtained from NMPC using the node-only and the CTCS formulations throughout the simulation. The legend reports the average discrepancy over the entire simulation.
The second plot shows the violation in the obstacle avoidance constraint, quantified by $1^\top |g(t,\xi(t),u(t))|_+$, throughout the simulation.}
\label{fig:tracking_error_3}
\end{figure}

We compare the proposed approach against the usual method of considering only the node point state and control inputs, in order to quantify the integral of the running cost and violation of the obstacle avoidance constraint. 
The output reference tracking and obstacle avoidance results are shown in Figures \ref{fig:sim_path_thick} and \ref{fig:tracking_error_3}.

Figure \ref{fig:sim_path_thick} indicates that in the node-only case, the obstacle avoidance constraint is not satisfied for the continuous-time trajectory (obtained via a single shooting propagation of the system dynamics with the NMPC control input), despite being satisfied at all node points. 
Conversely, the CTCS formulation ensures that the constraint is satisfied at all times along the trajectory.

Figure \ref{fig:tracking_error_3} demonstrates that both the node-only and CTCS formulations exhibit comparable performances in terms of output reference tracking. The CTCS formulation shows a slightly lower performance, but this difference is because the CTCS formulation actually satisfies the obstacle avoidance constraint.

\section{Conclusions}

\vspace{-0.15cm}
This paper presents the adaptation of a successive convexification framework to address the NMPC problem.
Reformulating the optimal control problem by embedding path constraints and the running cost into the dynamics function within the successive convexification framework offers the benefit of being able to operate on sparse discretization grids without encountering inter-sample constraint violations or a degradation in the accuracy of the running cost. The prox-linear method, a convergence-guaranteed SCP algorithm that requires only first-order information and supports warm-starting, is employed to solve the reformulated optimal control problem. 
The aforementioned attributes allow the algorithm to operate in an RTI-like mode while achieving a favorable solution in a few iterations and maintaining real-time capability.
The proposed framework is demonstrated on a tracking problem involving obstacle avoidance.


\section*{Acknowledgments}
\vspace{-0.15cm}
This research was supported by AFOSR grant FA9550-20-1-0053 and ONR grant N00014-20-1-2288; Government sponsorship is acknowledged.


\bibliography{ifacconf}

\end{document}